\documentclass[graybox]{svmult}
\spnewtheorem{theoremA}{Theorem}{\bfseries}{\itshape}

\usepackage{tikz}
\usepackage{type1cm}        
\usepackage{makeidx}         
\usepackage{graphicx}        
\usepackage{multicol}        
\usepackage[bottom]{footmisc}

\newenvironment{Proof}{\textbf{\textit{{Proof:}}}}{\hfill$\square$}

\usepackage{newtxtext}       %
\usepackage{newtxmath}       

\makeindex             

\newcommand{\PP}{{\mathbb P}}

\newcommand{\C}{{\mathbb C}}

\newcommand{\Z}{{\mathbb Z}}

\newcommand{\N}{{\mathbb N}}
\newcommand{\F}{{\mathbb F}}
\DeclareMathOperator\supp{supp}

\newcommand{\CP}{\C\PP}

\begin{document}

\title*{On the Holomorphic and Random Dynamics for some examples of higher rank Free Groups generated by H\'enon type maps\thanks{This is the author's accepted manuscript of the contribution. The Version of Record of this contribution will be published in [Analysis and PDE in Developing Countries: Proceedings of the ISAAC-ICMAM Conference, 2025]. The Version of Record will be made available online at [https://doi.org/[insert DOI]] upon publication.}}
\titlerunning{Examples of higher rank Free Groups generated by H\'enon type maps} 
\author{Andres Quintero Santander}
\institute{Andres Quintero Santander \at Indiana University Indianapolis,  420 University Boulevard Indianapolis, IN 46202 \email{aequinte@iu.edu}}
%
%
\maketitle

\date{\today}


\abstract*{Each chapter should be preceded by an abstract (no more than 200 words) that summarizes the content. The abstract will appear \textit{online} at \url{www.SpringerLink.com} and be available with unrestricted access. This allows unregistered users to read the abstract as a teaser for the complete chapter.
Please use the 'starred' version of the \texttt{abstract} command for typesetting the text of the online abstracts (cf. source file of this chapter template \texttt{abstract}) and include them with the source files of your manuscript. Use the plain \texttt{abstract} command if the abstract is also to appear in the printed version of the book.}

\abstract{We study the Holomorphic and Random Dynamics of some rank 2 free groups generated by two H\'enon type maps. For these simply constructed  examples we prove that the Fatou set is non-empty and that the stationary measures are supported on a compact set. With some further care this allows us to construct examples  having no stationary measures. These examples illustrate the types of phenomena that may arise when studying holomorphic group actions on non-compact manifolds. }
\section{Introduction}

Classical Holomorphic Dynamics typically considers iteration of a single mapping. In the literature, there are many examples  of mappings of non-compact complex manifolds. The first examples studied were polynomial maps from $\C$ to itself and on higher dimensions the first ones were the H\'enon maps on $\C^2$ see for example \cite{Hubbard1986}, \cite{bedford/Smillie} and \cite{MR1150591}. Actions of more complicated groups\footnote{For the remainder of the paper ''group action'' will refer to a group more complicated than $\Z$.} (or semigroups) of mappings lead to new interesting phenomena. 

Historically this has been done on compact manifolds/surfaces (see  e.g. \cite{CantatK3}). More recently, compactness of the underlying manifold has played an important role in the deep recent works of Cantat-Dujardin, see \cite{Cantat/Dujardin} and the references therein. On the other hand there are many natural examples of holomorphic group actions over non-compact manifolds, including the popular example of the Markov-type surfaces as it appears on \cite{cantat2009}, \cite{CantatMarkov}, \cite{cantatstationary}  and \cite{Rebelo/Roeder}. In such cases non-compactness leads to interesting situations e.g. presence of a non-empty Fatou set for most parameters (see Theorem E on \cite{Rebelo/Roeder} and also section 4.1 on \cite{cantat2009}) or absence of stationary measures (see the Main Theorem and Theorem 8.3 on \cite{cantatstationary}). On this manuscript we study some simple examples constructed using classical H\'enon maps of $\C^2$ illustrating both pheonomena.

\begin{definition}
    Let $P\in \C[y]$ with $\deg P\geq 2$ and $\delta\in\C\setminus\{0\}.$  A  generalized Hénon map is a map over $\C^2$ with the form: 

   \begin{equation*}
       h(x,y)=\left(y,P(y)-\delta x\right).
   \end{equation*}
\end{definition}




 Any affine conjugation of a finite composition of generalized H\'enon maps is called a H\'enon type map. In this work we'll consider the following group generated by H\'enon type maps: 
 
 \begin{equation}\label{Eq: G}
  G=\langle h_1,h_2\rangle\,\,\, \mbox{where}\,\, h_1(x,y)=\left(y,P(y)-\delta x\right)\,\, \mbox{and} \,\, h_2:=R_{\theta}^{-1}\circ h_1\circ R_{\theta}. 
 \end{equation}
 Here, $R_{\theta}:\C^2\to\C^2$ is the linear transformation induced by the matrix $\begin{pmatrix}
\cos\theta & -\sin\theta \\
\sin\theta & \cos\theta
\end{pmatrix}$ with $\theta\in(0,2\pi)\setminus\{\frac{\pi}{2},\frac{3\pi}{2},\pi\}$; this condition will imply that the group $G$ is not isomorphic to $\Z$, indeed, we will see later that $G\cong F_2$. We study the action of this group on $\C^2$ and the behavior of its orbits. We will prove that the Fatou set of $G$ is non-empty and explore the Random Dynamics and Ergodic Theory induced by the group $G$. We are interested in questions regarding sets invariant by $G$, orbit closures and also more quantitative questions about the statistics of random orbits. For the latter we need to describe the stationary measures. This concept is a common place for the study of random walks (see for example \cite{cantatstationary}). In particular, we give sufficient conditions for the absence of stationary measures to occur, in other words, conditions for the Random Dynamics to be trivial. This indicates a need for other methods of studying orbit closures for $G$.
\label{sec:1}
\section{H\'enon Maps, Complex Dynamics and Fatou-Julia Decomposition}
The dynamics of a single H\'enon map  (on $\C^2$) is described by the following filtration that appears on \cite{bedford-henon}:

\vspace{0.3 cm}

\tikzset{every picture/.style={line width=0.75pt}} 

\begin{tikzpicture}[x=0.75pt,y=0.75pt,yscale=-1,xscale=1]

\draw  (177,204.03) -- (449,204.03)(204.2,26.73) -- (204.2,223.73) (442,199.03) -- (449,204.03) -- (442,209.03) (199.2,33.73) -- (204.2,26.73) -- (209.2,33.73)  ;
\draw   (204.2,112.73) -- (318,112.73) -- (318,204.03) -- (204.2,204.03) -- cycle ;
\draw    (318,112.73) -- (420,22.73) ;
\draw  [dash pattern={on 0.84pt off 2.51pt}]  (368,149) .. controls (404.45,184.2) and (435.07,162.41) .. (429.29,123.52) ;
\draw [shift={(429,121.73)}, rotate = 80.07] [color={rgb, 255:red, 0; green, 0; blue, 0 }  ][line width=0.75]    (10.93,-3.29) .. controls (6.95,-1.4) and (3.31,-0.3) .. (0,0) .. controls (3.31,0.3) and (6.95,1.4) .. (10.93,3.29)   ;
\draw    (279,191.73) .. controls (301.77,167.97) and (290.24,132.45) .. (249.25,169.59) ;
\draw [shift={(248,170.73)}, rotate = 317.12] [color={rgb, 255:red, 0; green, 0; blue, 0 }  ][line width=0.75]    (10.93,-3.29) .. controls (6.95,-1.4) and (3.31,-0.3) .. (0,0) .. controls (3.31,0.3) and (6.95,1.4) .. (10.93,3.29)   ;
\draw    (238,156.73) -- (238,87.73) ;
\draw [shift={(238,85.73)}, rotate = 90] [color={rgb, 255:red, 0; green, 0; blue, 0 }  ][line width=0.75]    (10.93,-3.29) .. controls (6.95,-1.4) and (3.31,-0.3) .. (0,0) .. controls (3.31,0.3) and (6.95,1.4) .. (10.93,3.29)   ;
\draw    (354,165.73) -- (303,164.77) ;
\draw [shift={(301,164.73)}, rotate = 1.08] [color={rgb, 255:red, 0; green, 0; blue, 0 }  ][line width=0.75]    (10.93,-3.29) .. controls (6.95,-1.4) and (3.31,-0.3) .. (0,0) .. controls (3.31,0.3) and (6.95,1.4) .. (10.93,3.29)   ;
\draw    (373,103.73) -- (322.54,62.01) ;
\draw [shift={(321,60.73)}, rotate = 39.59] [color={rgb, 255:red, 0; green, 0; blue, 0 }  ][line width=0.75]    (10.93,-3.29) .. controls (6.95,-1.4) and (3.31,-0.3) .. (0,0) .. controls (3.31,0.3) and (6.95,1.4) .. (10.93,3.29)   ;
\draw    (289,97.73) -- (289,17.73) ;
\draw [shift={(289,15.73)}, rotate = 90] [color={rgb, 255:red, 0; green, 0; blue, 0 }  ][line width=0.75]    (10.93,-3.29) .. controls (6.95,-1.4) and (3.31,-0.3) .. (0,0) .. controls (3.31,0.3) and (6.95,1.4) .. (10.93,3.29)   ;

\draw (176,18.4) node [anchor=north west][inner sep=0.75pt]    {$|y|$};
\draw (455,199.4) node [anchor=north west][inner sep=0.75pt]    {$|x|$};
\draw (311,209.4) node [anchor=north west][inner sep=0.75pt]    {$R$};
\draw (261,11.4) node [anchor=north west][inner sep=0.75pt]    {$\infty $};
\draw (231,43.4) node [anchor=north west][inner sep=0.75pt]    {$V^{-}$};
\draw (403,96.4) node [anchor=north west][inner sep=0.75pt]    {$V^{+}$};
\draw (396,4.4) node [anchor=north west][inner sep=0.75pt]    {$|y|=|x|$};
\draw (218,158.4) node [anchor=north west][inner sep=0.75pt]    {$V$};
\draw (183,105.4) node [anchor=north west][inner sep=0.75pt]    {$R$};

\end{tikzpicture}

Where:
\begin{enumerate}
    \item $h(V\cup V^-)\subset  V\cup V^-;$
    \item $h(V^-)\subset V^-;$
    \item For all $(x,y)\in V^-$  we have $\displaystyle \lim_{n \to \infty}\left\| h^{n}(x,y)\right\|=\infty ;$ and
    \item The dashed arrow represents that the forward orbit of a point $(x,y)\in V^+$ under $h$ can only remain in $V^+$ a finite amount of time and after that this orbit has to go either to $V$ or enter the $V^-$ region and escape to infinity.
\end{enumerate}

A classical tool in Complex Dynamics is the Fatou set (and its counter-part the Julia set).  We define normality and the Fatou set next (cf. \cite{Zalcman}):

\begin{definition}
    A family of maps $\mathcal{A}=\{f_\alpha:X\to Y| \alpha\in I\}$ is called normal if for every sequence $(f_n)_{n\in\N}\subset \mathcal{A}$ either there is a subsequence that converges uniformly on compact subsets to some function $f:X\to Y$ or properly diverges to infinity. The latter is only needed if $Y$ is non-compact.
\end{definition}
\begin{definition}
    The Fatou set $\mathcal{F}(\mathcal{A})$ is the largest open subset of the domain $X$ such that the family $\mathcal{A}=\{f_\alpha:X\to Y| \alpha\in I\}$ is normal. The Julia set $\mathcal{J}(\mathcal{A})$ of $\mathcal{A}$ is defined as the complement of the Fatou set.
\end{definition}
The study of the Fatou set  has been typically associated with the set of all iterates of a single map and this has been studied a lot in the literature, but as we saw in the definition this can also be extended for more general families of maps like groups or semigroups of mappings or even more generally.
The main idea is that the Julia set represents the region  with chaotic behavior under the iterates of the (semi)group action. In the Julia set, orbit behavior is difficult to predict and sensitive to initial conditions. Conversely, the Fatou set is where the dynamics are more stable. 
\begin{theoremA}
Let $G$ be defined as in \eqref{Eq: G}. Then there is some unbounded open set $\emptyset\neq U\subset \mathcal{F}(G)$ such that the orbit of every $(x,y)\in U$ escapes to infinity under every sequence $(\omega_n)_{n\in\N}\in\{h_1,h_2,h_1^{-1},h_2^{-1}\}^\N$ such that $\mbox{length}(S_k(\omega)) \underset{k\to\infty}{\to}\infty.$ Here $S_k(\omega):=\omega_k\circ\omega_{k-1}\circ\dots \circ\omega_1$ and we consider its length after reductions. 
\end{theoremA}
\begin{remark}
 Note that Theorem A implies that there is no $p\in \C^2$ whose orbit is dense in $\C^2$. It would be interesting to know if there is a point $p$ in the Julia set $\mathcal{J}(G)$ whose orbit is dense in $\mathcal{J}(G)$ (Cf. Theorem A of \cite{Rebelo/Roeder}) 
\end{remark}

Since Hénon maps are invertible, in principle, we can study the group generated by two arbitrary Hénon maps and explore the dynamics of its group action on $\C^2$. We have reduced the scope of the question by assuming the second generator $h_2$ is the conjugate of $h_1$ via a rotation. Furthermore, we require that this rotation is by $\theta\in(0,2\pi)\setminus \{\frac{\pi}{2},\frac{3\pi}{2},\pi\}$ so that it separates the indeterminacies of $h_1, h_2$, and their inverses as maps of $\CP^2$.  More precisely, let $h_1(x,y)=\left(y,P(y)-\delta x\right)$ with $\deg(P)=d\geq 2$ be a generalized H\'enon map, then if $\tilde{P}(x_1,x_2)=x_2^dP(x_1/x_2)$, we have that $h_1$ can be extended as a map of the projective plane $\CP^2$ (in homogeneous coordinates) as:
\begin{equation}\label{eq: Henon map 1} 
    h_1[x_0\colon x_1\colon x_2]=[x_1x_2^{d-1}\colon \tilde{P}(x_1,x_2)-\delta x_0x_2^{d-1}\colon x_2^d].
\end{equation}

We have gained compactness for the phase space but the price is that now $h_1$ is not defined everywhere. Let $L_\infty=\{x_2=0\}$ be the line at infinity, notice, that the point $p=[1:0:0]\in L_\infty$ is the only indeterminate point of this birational map and that the line at infinity $L_\infty\setminus\{p\}$ is mapped by $h_1$ to the fixed point $q=[0:1:0]$. Similarly, $q$ is the  only indetermination of $h_1^{-1}$ and $h_1^{-1}$ contracts $L_\infty\setminus\{q\}$ to $p$. As $R_\theta$ is a rotation, we can compute its homogeneous coordinates and notice that $h_2$ has also a super-attracting fixed point $\tilde{q}$ on the line at infinity and a unique indetermination $\tilde{p}$ there also and analogously as before the behavior of its inverse its the opposite. So, we pick $\theta\in(0,2\pi)\setminus \{\frac{\pi}{2},\frac{3\pi}{2},\pi\}$ so that the indeterminations $\tilde{p}$ of $h_2:=R_\theta^{-1}\circ h_1\circ R_\theta$ and $\tilde{q}$ of $h_2^{-1}$ are such that the set $\{p,q,\tilde{p},\tilde{q}\}$ has 4 distinct elements. As we said before we are interested in showing that the Fatou set for the group $G=\langle h_1,h_2\rangle$ on this context is non-empty. 

Our strategy is the following: we claim the existence of an unbounded component of the Fatou set $\mathcal{F}(G)$ that contains a neighborhood $U$ of $L_{\infty}\setminus\{p,q,\tilde{p},\tilde{q}\}$ such that for every $(x,y)\in U\cap \C^2$ and for every $\omega=(\omega_{n})_{n\in\N}\subset \{h_1,h_2,h_{1}^{-1},h_{2}^{-1}\}^\N$ such that $\mbox{length}(S_k(\omega)) \underset{k\to\infty}{\to}\infty$ we have $\displaystyle\lim_{k\to\infty}\|S_k(\omega)\cdot (x,y)\|=\infty$. The following is an illustration of $\CP^2$ together with the subset of the Fatou set $U$: 

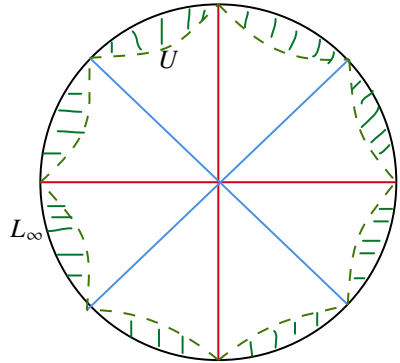
\begin{center}

\tikzset{every picture/.style={line width=0.75pt}} 

\begin{tikzpicture}[x=0.75pt,y=0.75pt,yscale=-0.7,xscale=0.7]

\draw [color={rgb, 255:red, 208; green, 2; blue, 27 }  ,draw opacity=1 ]   (307.76,23) -- (307.76,279.48) ;
\draw   (179.52,151.24) .. controls (179.52,80.42) and (236.93,23) .. (307.76,23) .. controls (378.58,23) and (436,80.42) .. (436,151.24) .. controls (436,222.07) and (378.58,279.48) .. (307.76,279.48) .. controls (236.93,279.48) and (179.52,222.07) .. (179.52,151.24) -- cycle ;
\draw [color={rgb, 255:red, 208; green, 2; blue, 27 }  ,draw opacity=1 ]   (179.52,151.24) -- (436,151.24) ;
\draw [color={rgb, 255:red, 74; green, 144; blue, 226 }  ,draw opacity=1 ]   (402,62.48) -- (213,243.48) ;
\draw [color={rgb, 255:red, 74; green, 144; blue, 226 }  ,draw opacity=1 ]   (215,61.48) -- (221.11,67.33) -- (401,239.48) ;
\draw [color={rgb, 255:red, 65; green, 117; blue, 5 }  ,draw opacity=1 ] [dash pattern={on 4.5pt off 4.5pt}]  (215,61.48) .. controls (277,58.48) and (274,56.48) .. (307.76,23) ;
\draw [color={rgb, 255:red, 65; green, 117; blue, 5 }  ,draw opacity=1 ] [dash pattern={on 4.5pt off 4.5pt}]  (307.76,23) .. controls (334,49.48) and (354,66.48) .. (402,62.48) ;
\draw [color={rgb, 255:red, 65; green, 117; blue, 5 }  ,draw opacity=1 ] [dash pattern={on 4.5pt off 4.5pt}]  (179.52,151.24) .. controls (222,112.48) and (214,103.48) .. (215,61.48) ;
\draw [color={rgb, 255:red, 65; green, 117; blue, 5 }  ,draw opacity=1 ] [dash pattern={on 4.5pt off 4.5pt}]  (213,243.48) .. controls (217,195.48) and (215,189.48) .. (179.52,151.24) ;
\draw [color={rgb, 255:red, 65; green, 117; blue, 5 }  ,draw opacity=1 ] [dash pattern={on 4.5pt off 4.5pt}]  (307.76,279.48) .. controls (285,260.48) and (254,242.48) .. (213,243.48) ;
\draw [color={rgb, 255:red, 65; green, 117; blue, 5 }  ,draw opacity=1 ] [dash pattern={on 4.5pt off 4.5pt}]  (307.76,279.48) .. controls (324.9,266.01) and (338.54,258.34) .. (350.58,253.44) .. controls (368.9,245.98) and (383.49,244.92) .. (401,239.48) ;
\draw [color={rgb, 255:red, 65; green, 117; blue, 5 }  ,draw opacity=1 ] [dash pattern={on 4.5pt off 4.5pt}]  (402,62.48) .. controls (406,112.48) and (407,117.48) .. (436,151.24) ;
\draw [color={rgb, 255:red, 65; green, 117; blue, 5 }  ,draw opacity=1 ] [dash pattern={on 4.5pt off 4.5pt}]  (401,239.48) .. controls (404,187.48) and (410,183.48) .. (436.48,149.73) ;
\draw  [color={rgb, 255:red, 5; green, 117; blue, 43 }  ,draw opacity=1 ][line width=0.75] [line join = round][line cap = round] (231,50.48) .. controls (231,52.82) and (231,55.15) .. (231,57.48) ;
\draw  [color={rgb, 255:red, 5; green, 117; blue, 43 }  ,draw opacity=1 ][line width=0.75] [line join = round][line cap = round] (246,42.48) .. controls (250.48,46.96) and (247.95,51.44) .. (252,55.48) ;
\draw  [color={rgb, 255:red, 5; green, 117; blue, 43 }  ,draw opacity=1 ][line width=0.75] [line join = round][line cap = round] (267,33.48) .. controls (267,39.48) and (267,45.48) .. (267,51.48) ;
\draw  [color={rgb, 255:red, 5; green, 117; blue, 43 }  ,draw opacity=1 ][line width=0.75] [line join = round][line cap = round] (287,26.48) .. controls (287,32.17) and (286,35.08) .. (286,40.48) ;
\draw  [color={rgb, 255:red, 5; green, 117; blue, 43 }  ,draw opacity=1 ][line width=0.75] [line join = round][line cap = round] (298,27.48) .. controls (296.77,29.95) and (297.23,32.02) .. (296,34.48) ;
\draw  [color={rgb, 255:red, 5; green, 117; blue, 43 }  ,draw opacity=1 ][line width=0.75] [line join = round][line cap = round] (200,87.48) .. controls (203,87.48) and (206,87.48) .. (209,87.48) ;
\draw  [color={rgb, 255:red, 5; green, 117; blue, 43 }  ,draw opacity=1 ][line width=0.75] [line join = round][line cap = round] (193,99.48) .. controls (198.34,99.48) and (203.66,100.48) .. (209,100.48) ;
\draw  [color={rgb, 255:red, 5; green, 117; blue, 43 }  ,draw opacity=1 ][line width=0.75] [line join = round][line cap = round] (192,114.48) .. controls (198.83,114.48) and (204.03,115.48) .. (211,115.48) ;
\draw  [color={rgb, 255:red, 5; green, 117; blue, 43 }  ,draw opacity=1 ][line width=0.75] [line join = round][line cap = round] (182,130.48) .. controls (186,130.48) and (190,130.48) .. (194,130.48) ;
\draw  [color={rgb, 255:red, 5; green, 117; blue, 43 }  ,draw opacity=1 ][line width=0.75] [line join = round][line cap = round] (185,186.48) .. controls (191.91,186.48) and (194.19,184.48) .. (202,184.48) ;
\draw  [color={rgb, 255:red, 5; green, 117; blue, 43 }  ,draw opacity=1 ][line width=0.75] [line join = round][line cap = round] (188,176.48) .. controls (192.67,176.48) and (197.33,176.48) .. (202,176.48) ;
\draw  [color={rgb, 255:red, 5; green, 117; blue, 43 }  ,draw opacity=1 ][line width=0.75] [line join = round][line cap = round] (192,203.48) .. controls (198,203.48) and (204,203.48) .. (210,203.48) ;
\draw  [color={rgb, 255:red, 5; green, 117; blue, 43 }  ,draw opacity=1 ][line width=0.75] [line join = round][line cap = round] (185,168.48) .. controls (189,168.48) and (193,168.48) .. (197,168.48) ;
\draw  [color={rgb, 255:red, 5; green, 117; blue, 43 }  ,draw opacity=1 ][line width=0.75] [line join = round][line cap = round] (199,218.48) .. controls (202.33,218.48) and (205.67,218.48) .. (209,218.48) ;
\draw  [color={rgb, 255:red, 5; green, 117; blue, 43 }  ,draw opacity=1 ][line width=0.75] [line join = round][line cap = round] (209,218.48) .. controls (209,218.48) and (209,218.48) .. (209,218.48) ;
\draw  [color={rgb, 255:red, 5; green, 117; blue, 43 }  ,draw opacity=1 ][line width=0.75] [line join = round][line cap = round] (245,251.48) .. controls (245,255.15) and (245,258.82) .. (245,262.48) ;
\draw  [color={rgb, 255:red, 5; green, 117; blue, 43 }  ,draw opacity=1 ][line width=0.75] [line join = round][line cap = round] (265,257.48) .. controls (265,261.48) and (265,265.48) .. (265,269.48) ;
\draw  [color={rgb, 255:red, 5; green, 117; blue, 43 }  ,draw opacity=1 ][line width=0.75] [line join = round][line cap = round] (282,262.48) .. controls (282,266.48) and (282,270.48) .. (282,274.48) ;
\draw  [color={rgb, 255:red, 5; green, 117; blue, 43 }  ,draw opacity=1 ][line width=0.75] [line join = round][line cap = round] (331,266.48) .. controls (331,269.83) and (332,273.13) .. (332,276.48) ;
\draw  [color={rgb, 255:red, 5; green, 117; blue, 43 }  ,draw opacity=1 ][line width=0.75] [line join = round][line cap = round] (348,258.48) .. controls (350.84,261.32) and (349,266.47) .. (349,270.48) ;
\draw  [color={rgb, 255:red, 5; green, 117; blue, 43 }  ,draw opacity=1 ][line width=0.75] [line join = round][line cap = round] (363,252.48) .. controls (363,255.15) and (363,257.82) .. (363,260.48) ;
\draw  [color={rgb, 255:red, 5; green, 117; blue, 43 }  ,draw opacity=1 ][line width=0.75] [line join = round][line cap = round] (379,245.48) .. controls (379,248.48) and (379,251.48) .. (379,254.48) ;
\draw  [color={rgb, 255:red, 5; green, 117; blue, 43 }  ,draw opacity=1 ][line width=0.75] [line join = round][line cap = round] (406,220.48) .. controls (408.33,220.48) and (410.67,220.48) .. (413,220.48) ;
\draw  [color={rgb, 255:red, 5; green, 117; blue, 43 }  ,draw opacity=1 ][line width=0.75] [line join = round][line cap = round] (410,202.48) .. controls (413.33,202.48) and (416.67,202.48) .. (420,202.48) ;
\draw  [color={rgb, 255:red, 5; green, 117; blue, 43 }  ,draw opacity=1 ][line width=0.75] [line join = round][line cap = round] (412,193.48) .. controls (416,193.48) and (420,193.48) .. (424,193.48) ;
\draw  [color={rgb, 255:red, 5; green, 117; blue, 43 }  ,draw opacity=1 ][line width=0.75] [line join = round][line cap = round] (420,179.48) .. controls (423.33,179.48) and (426.67,179.48) .. (430,179.48) ;
\draw  [color={rgb, 255:red, 5; green, 117; blue, 43 }  ,draw opacity=1 ][line width=0.75] [line join = round][line cap = round] (424,170.48) .. controls (427,170.48) and (430,170.48) .. (433,170.48) ;
\draw  [color={rgb, 255:red, 5; green, 117; blue, 43 }  ,draw opacity=1 ][line width=0.75] [line join = round][line cap = round] (410,78.48) .. controls (408,80.48) and (406,82.48) .. (404,84.48) ;
\draw  [color={rgb, 255:red, 5; green, 117; blue, 43 }  ,draw opacity=1 ][line width=0.75] [line join = round][line cap = round] (419,85.48) .. controls (419,88.28) and (411,89.78) .. (411,94.48) ;
\draw  [color={rgb, 255:red, 5; green, 117; blue, 43 }  ,draw opacity=1 ][line width=0.75] [line join = round][line cap = round] (423,97.48) .. controls (423,101.38) and (413.87,108.61) .. (411,111.48) ;
\draw  [color={rgb, 255:red, 5; green, 117; blue, 43 }  ,draw opacity=1 ][line width=0.75] [line join = round][line cap = round] (428,113.48) .. controls (423.46,118.03) and (422.54,117.94) .. (418,122.48) ;
\draw  [color={rgb, 255:red, 5; green, 117; blue, 43 }  ,draw opacity=1 ][line width=0.75] [line join = round][line cap = round] (431,124.48) .. controls (429.53,125.95) and (425.77,130.48) .. (424,130.48) ;
\draw  [color={rgb, 255:red, 5; green, 117; blue, 43 }  ,draw opacity=1 ][line width=0.75] [line join = round][line cap = round] (329,28.48) .. controls (329,31.57) and (327,34.24) .. (327,39.48) ;
\draw  [color={rgb, 255:red, 5; green, 117; blue, 43 }  ,draw opacity=1 ][line width=0.75] [line join = round][line cap = round] (349,33.48) .. controls (349,39.08) and (347.32,41.22) .. (346,46.48) ;
\draw  [color={rgb, 255:red, 5; green, 117; blue, 43 }  ,draw opacity=1 ][line width=0.75] [line join = round][line cap = round] (364,42.48) .. controls (364,46.88) and (361.21,51.07) .. (359,55.48) ;
\draw  [color={rgb, 255:red, 5; green, 117; blue, 43 }  ,draw opacity=1 ][line width=0.75] [line join = round][line cap = round] (381,48.48) .. controls (381,50.83) and (379.26,60.48) .. (376,60.48) ;
\draw  [color={rgb, 255:red, 5; green, 117; blue, 43 }  ,draw opacity=1 ][line width=0.75] [line join = round][line cap = round] (391,57.48) .. controls (391,60.69) and (389.03,63.45) .. (387,65.48) ;

\draw (263,54.4) node [anchor=north west][inner sep=0.75pt]    {$U$};
\draw (155,176.4) node [anchor=north west][inner sep=0.75pt]    {$L_{\infty }$};

\end{tikzpicture}
\end{center}

\label{subsec:2}

 Here the boundary of the circle is the line at infinity, the set in green is our subset $U$ of the Fatou set, the end-points of the red rays  represent the indeterminacies of $h_1$ and $h_1^{-1}$ and the end-points of the blue rays the indeterminacies of $h_2$ and $h_2^{-1}.$ The picture was created by assuming that we rotated by $R_{\frac{\pi}{4}}$ but the reader can notice that our reasoning works for any rotation that separates the indeterminacies of $h_1$ and $h_2$ and their inverses.

First, we highlight that proving the existence of a Fatou set for the Markov-type examples required a more delicate treatment and involved ideas coming from Teichmüller Theory; see for example the proof and commentaries of Theorem E in \cite{Rebelo/Roeder}. To prove Theorem A however, we only needed to use topological ideas.

\begin{Proof}
First consider a set $W_1\subset V^{-}_{h_2}$ (i.e. escaping to infinity by iteration of $h_2$) such that $\|h_1(x,y)\|>2\|(x,y)\|$ and $\|h_1^{-1}(x,y)\|>2\|(x,y)\|$ for all $(x,y)\in W_1$ (also for $h_2$ by taking it far out enough). Said differently, $W_1$ is an open set where every point increases more than twice its norm under 3 out of the 4 maps. This set clearly exists, just pick a truncated cone around the attracting direction of $h_2$ and require for that cone to be far enough from the indeterminacies of $h_1$ and $h_1^{-1}$. Indeed, for this we can use the fact that $\deg(P)\geq 2$ to make its norm increase under  both $h_1$ and $h_1^{-1}$. 
\vspace{0.2 cm}

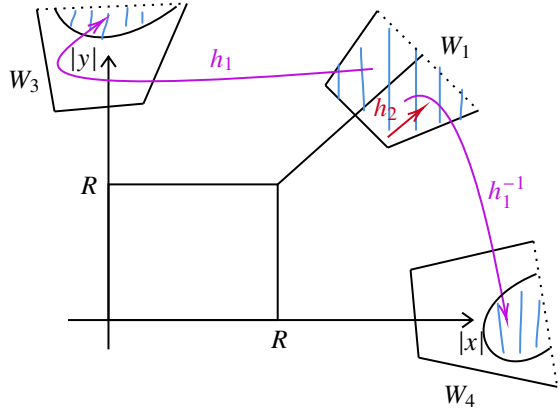
\tikzset{every picture/.style={line width=0.75pt}} 

\begin{tikzpicture}[x=0.75pt,y=0.75pt,yscale=-0.75,xscale=0.75]

\draw  (138,236.03) -- (410,236.03)(165.2,58.73) -- (165.2,255.73) (403,231.03) -- (410,236.03) -- (403,241.03) (160.2,65.73) -- (165.2,58.73) -- (170.2,65.73)  ;
\draw   (165.2,144.73) -- (279,144.73) -- (279,236.03) -- (165.2,236.03) -- cycle ;
\draw    (279,144.73) -- (376,58) ;
\draw    (311,78) -- (353,120) ;
\draw    (311,78) -- (346,28) ;
\draw    (353,120) -- (414,95) ;
\draw    (368.34,202.16) -- (373.79,261.3) ;
\draw    (368.34,202.16) -- (451,183) ;
\draw    (373.79,261.3) -- (467,281) ;
\draw    (128.93,95.66) -- (188.72,90.36) ;
\draw    (128.93,95.66) -- (117,27) ;
\draw    (188.72,90.36) -- (218,23) ;
\draw  [color={rgb, 255:red, 74; green, 144; blue, 226 }  ,draw opacity=1 ][line width=0.75] [line join = round][line cap = round] (335,53) .. controls (335,67.12) and (337,78.7) .. (337,95) ;
\draw  [color={rgb, 255:red, 74; green, 144; blue, 226 }  ,draw opacity=1 ][line width=0.75] [line join = round][line cap = round] (354,40) .. controls (354,62.67) and (354,85.33) .. (354,108) ;
\draw  [color={rgb, 255:red, 74; green, 144; blue, 226 }  ,draw opacity=1 ][line width=0.75] [line join = round][line cap = round] (372,53) .. controls (372,71.33) and (372,89.67) .. (372,108) ;
\draw  [color={rgb, 255:red, 74; green, 144; blue, 226 }  ,draw opacity=1 ][line width=0.75] [line join = round][line cap = round] (388,66) .. controls (388,78.33) and (388,90.67) .. (388,103) ;
\draw  [color={rgb, 255:red, 74; green, 144; blue, 226 }  ,draw opacity=1 ][line width=0.75] [line join = round][line cap = round] (403,81) .. controls (403,90.13) and (402,93.6) .. (402,101) ;
\draw  [color={rgb, 255:red, 74; green, 144; blue, 226 }  ,draw opacity=1 ][line width=0.75] [line join = round][line cap = round] (320,64) .. controls (320,71.67) and (320,79.33) .. (320,87) ;
\draw [color={rgb, 255:red, 189; green, 16; blue, 224 }  ,draw opacity=1 ]   (364,89) .. controls (403.2,59.6) and (426.07,201.15) .. (432.62,234.12) ;
\draw [shift={(433,236)}, rotate = 258.31] [color={rgb, 255:red, 189; green, 16; blue, 224 }  ,draw opacity=1 ][line width=0.75]    (10.93,-3.29) .. controls (6.95,-1.4) and (3.31,-0.3) .. (0,0) .. controls (3.31,0.3) and (6.95,1.4) .. (10.93,3.29)   ;
\draw    (453,205) .. controls (390,231) and (422,285) .. (462,255) ;
\draw  [dash pattern={on 0.84pt off 2.51pt}]  (451,183) -- (467,281) ;
\draw  [color={rgb, 255:red, 74; green, 144; blue, 226 }  ,draw opacity=1 ][line width=0.75] [line join = round][line cap = round] (427,225) .. controls (427,226.31) and (430.39,257.39) .. (431,258) ;
\draw  [color={rgb, 255:red, 74; green, 144; blue, 226 }  ,draw opacity=1 ][line width=0.75] [line join = round][line cap = round] (443,217) .. controls (443,231.85) and (442,242.88) .. (442,258) ;
\draw  [color={rgb, 255:red, 74; green, 144; blue, 226 }  ,draw opacity=1 ][line width=0.75] [line join = round][line cap = round] (455,223) .. controls (455,234.59) and (453,243.16) .. (453,254) ;
\draw    (129,26) .. controls (139,50) and (171,55) .. (211,25) ;
\draw  [dash pattern={on 0.84pt off 2.51pt}]  (117,27) -- (159.05,25.33) -- (218,23) ;
\draw  [dash pattern={on 0.84pt off 2.51pt}]  (346,28) -- (414,95) ;
\draw  [color={rgb, 255:red, 74; green, 144; blue, 226 }  ,draw opacity=1 ][line width=0.75] [line join = round][line cap = round] (143,28) .. controls (143,31.51) and (140,34.13) .. (140,38) ;
\draw  [color={rgb, 255:red, 74; green, 144; blue, 226 }  ,draw opacity=1 ][line width=0.75] [line join = round][line cap = round] (154,28) .. controls (154,32.95) and (153,35.62) .. (153,41) ;
\draw  [color={rgb, 255:red, 74; green, 144; blue, 226 }  ,draw opacity=1 ][line width=0.75] [line join = round][line cap = round] (165,30) .. controls (165,36.54) and (161,41.33) .. (161,47) ;
\draw  [color={rgb, 255:red, 74; green, 144; blue, 226 }  ,draw opacity=1 ][line width=0.75] [line join = round][line cap = round] (178,26) .. controls (178,30.01) and (177,33.99) .. (177,38) ;
\draw  [color={rgb, 255:red, 74; green, 144; blue, 226 }  ,draw opacity=1 ][line width=0.75] [line join = round][line cap = round] (188,29) .. controls (188,31.33) and (188,33.67) .. (188,36) ;
\draw [color={rgb, 255:red, 189; green, 16; blue, 224 }  ,draw opacity=1 ]   (343,67) .. controls (77.68,88.78) and (121.09,64.5) .. (160.8,34.9) ;
\draw [shift={(162,34)}, rotate = 143.13] [color={rgb, 255:red, 189; green, 16; blue, 224 }  ,draw opacity=1 ][line width=0.75]    (10.93,-3.29) .. controls (6.95,-1.4) and (3.31,-0.3) .. (0,0) .. controls (3.31,0.3) and (6.95,1.4) .. (10.93,3.29)   ;
\draw [color={rgb, 255:red, 208; green, 2; blue, 27 }  ,draw opacity=1 ]   (353,113) -- (375.46,94.28) ;
\draw [shift={(377,93)}, rotate = 140.19] [color={rgb, 255:red, 208; green, 2; blue, 27 }  ,draw opacity=1 ][line width=0.75]    (10.93,-3.29) .. controls (6.95,-1.4) and (3.31,-0.3) .. (0,0) .. controls (3.31,0.3) and (6.95,1.4) .. (10.93,3.29)   ;

\draw (137,50.4) node [anchor=north west][inner sep=0.75pt]    {$|y|$};
\draw (398,241.4) node [anchor=north west][inner sep=0.75pt]    {$|x|$};
\draw (272,241.4) node [anchor=north west][inner sep=0.75pt]    {$R$};
\draw (144,137.4) node [anchor=north west][inner sep=0.75pt]    {$R$};
\draw (386,43.4) node [anchor=north west][inner sep=0.75pt]    {$W_{1}$};
\draw (418,137.4) node [anchor=north west][inner sep=0.75pt]  [color={rgb, 255:red, 189; green, 16; blue, 224 }  ,opacity=1 ]  {$h_{1}^{-1}$};
\draw (231,51.4) node [anchor=north west][inner sep=0.75pt]  [color={rgb, 255:red, 189; green, 16; blue, 224 }  ,opacity=1 ]  {$h_{1}$};
\draw (341,84.4) node [anchor=north west][inner sep=0.75pt]  [color={rgb, 255:red, 208; green, 2; blue, 27 }  ,opacity=1 ]  {$h_{2}$};
\draw (389,278.4) node [anchor=north west][inner sep=0.75pt]    {$W_{4}$};
\draw (96,66.4) node [anchor=north west][inner sep=0.75pt]    {$W_{3}$};

\end{tikzpicture}

In an analogous way, there are truncated cones $W_2, W_3$ and $W_4$ such that points in there increase in norm under application of $\{h_2^{-1}, h_1,h_1^{-1}\}$, $\{h_1, h_2,h_2^{-1}\}$ and $\{h_1^{-1}, h_2,h_2^{-1}\}$ respectively. Without losing generality, due to the escaping behavior of H\'enon maps we can consider these sets to be periodic in the following sense: if $\{a,b,b^{-1}\}$ are the three maps that double the norm for points in $W_i$ (where $a$ gives the main attracting direction) then $b(W_i)\subset W_j$ where $W_j$ increases length under $\{b,a, a^{-1}\}$, analogously for $b^{-1}$. Such truncated cones are as in the illustration above. 
In particular, If $(x,y)\in W_1$ and $\omega_1\in\{h_2,h_1,h_1^{-1}\}$ and there are no trivial reductions in the composition $S_n(\omega)$ we have:
\begin{equation}\label{Eq: increasing W1}
    \|S_n(\omega)\cdot (x,y)\|>2^n\|(x,y)\|.
\end{equation}

Analogous results follow for each $W_i$. Now, consider the map $h_1$ which has $q$ as a super-attracting fixed point and $p$ is its indeterminate point. Then as $\CP^2$ is compact and  $h_1$ continuous, we can cut out a neighborhood $A$ of the indetermination $p$ and use uniform continuity out of it to get a neighborhood $B$ of the part of the line at infinity out of $A$ such that $h_1(B)\subset W_3$. Improve this set $B$ so that includes a bigger portion of $L_{\infty}$ by cutting out smaller neighborhoods of $p$ and repeating the argument. This way we can obtain by exhaustion an open set $V_1$ containing $L_\infty\setminus\{p\}$ such that $h_1(V_1)\subset W_3$, where $W_3$ is the truncated cone around $q$.

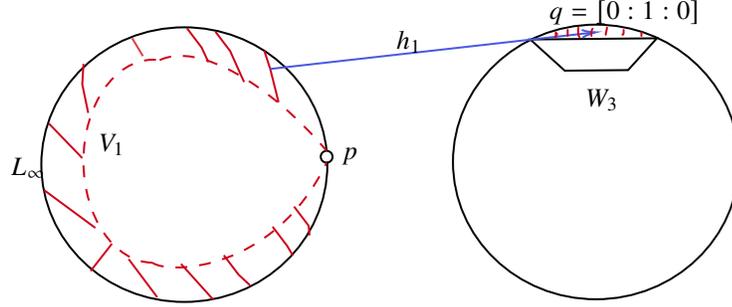
\tikzset{every picture/.style={line width=0.75pt}} 

\begin{tikzpicture}[x=0.75pt,y=0.75pt,yscale=-0.6,xscale=0.6]

\draw   (267.43,150.12) .. controls (267.43,147.38) and (269.65,145.15) .. (272.4,145.15) .. controls (275.15,145.15) and (277.38,147.38) .. (277.38,150.12) .. controls (277.38,152.87) and (275.15,155.1) .. (272.4,155.1) .. controls (269.65,155.1) and (267.43,152.87) .. (267.43,150.12) -- cycle ;
\draw  [draw opacity=0] (272.97,154.23) .. controls (272.99,155.07) and (273,155.91) .. (273,156.75) .. controls (273,220.54) and (219.18,272.25) .. (152.8,272.25) .. controls (86.41,272.25) and (32.6,220.54) .. (32.6,156.75) .. controls (32.6,92.96) and (86.41,41.25) .. (152.8,41.25) .. controls (215.11,41.25) and (266.35,86.81) .. (272.4,145.15) -- (152.8,156.75) -- cycle ; \draw   (272.97,154.23) .. controls (272.99,155.07) and (273,155.91) .. (273,156.75) .. controls (273,220.54) and (219.18,272.25) .. (152.8,272.25) .. controls (86.41,272.25) and (32.6,220.54) .. (32.6,156.75) .. controls (32.6,92.96) and (86.41,41.25) .. (152.8,41.25) .. controls (215.11,41.25) and (266.35,86.81) .. (272.4,145.15) ;  
\draw [color={rgb, 255:red, 208; green, 2; blue, 27 }  ,draw opacity=1 ] [dash pattern={on 4.5pt off 4.5pt}]  (132.07,68.87) .. controls (165.07,53.1) and (226,90.25) .. (272.4,145.15) ;
\draw [color={rgb, 255:red, 208; green, 2; blue, 27 }  ,draw opacity=1 ] [dash pattern={on 4.5pt off 4.5pt}]  (126.02,239.25) .. controls (139.02,249.25) and (228,247.25) .. (272.4,155.1) ;
\draw  [draw opacity=0][dash pattern={on 4.5pt off 4.5pt}] (126.04,239.25) .. controls (93.94,238.89) and (68,200.8) .. (68,153.87) .. controls (68,106.71) and (94.19,68.48) .. (126.5,68.48) .. controls (128.38,68.48) and (130.24,68.61) .. (132.07,68.87) -- (126.5,153.87) -- cycle ; \draw  [color={rgb, 255:red, 208; green, 2; blue, 27 }  ,draw opacity=1 ][dash pattern={on 4.5pt off 4.5pt}] (126.04,239.25) .. controls (93.94,238.89) and (68,200.8) .. (68,153.87) .. controls (68,106.71) and (94.19,68.48) .. (126.5,68.48) .. controls (128.38,68.48) and (130.24,68.61) .. (132.07,68.87) ;  
\draw [color={rgb, 255:red, 208; green, 2; blue, 27 }  ,draw opacity=1 ]   (35,178.25) -- (78,210) ;
\draw [color={rgb, 255:red, 208; green, 2; blue, 27 }  ,draw opacity=1 ]   (179,45) -- (191.32,61.25) -- (203,87.25) ;
\draw [color={rgb, 255:red, 208; green, 2; blue, 27 }  ,draw opacity=1 ]   (220,61) -- (225,76.25) -- (232,104.25) ;
\draw [color={rgb, 255:red, 208; green, 2; blue, 27 }  ,draw opacity=1 ]   (155,42) -- (167.32,58.25) -- (173,70.48) ;
\draw [color={rgb, 255:red, 208; green, 2; blue, 27 }  ,draw opacity=1 ]   (186,233) -- (198.32,249.25) -- (208,258.25) ;
\draw [color={rgb, 255:red, 208; green, 2; blue, 27 }  ,draw opacity=1 ]   (227,209) -- (239.32,225.25) -- (244,233.25) ;
\draw [color={rgb, 255:red, 208; green, 2; blue, 27 }  ,draw opacity=1 ]   (245,192) -- (257.32,208.25) -- (258,216.25) ;
\draw [color={rgb, 255:red, 208; green, 2; blue, 27 }  ,draw opacity=1 ]   (68,72) -- (67,89.25) -- (72,113.25) ;
\draw [color={rgb, 255:red, 208; green, 2; blue, 27 }  ,draw opacity=1 ]   (39,122.25) -- (67,151.25) ;
\draw [color={rgb, 255:red, 208; green, 2; blue, 27 }  ,draw opacity=1 ]   (106,236.25) -- (130,268.25) ;
\draw [color={rgb, 255:red, 208; green, 2; blue, 27 }  ,draw opacity=1 ]   (150,240.25) -- (177,270.25) ;
\draw [color={rgb, 255:red, 208; green, 2; blue, 27 }  ,draw opacity=1 ]   (92,223.25) -- (74,246.25) ;
\draw [color={rgb, 255:red, 218; green, 66; blue, 66 }  ,draw opacity=1 ]   (109,50.48) -- (118,65.48) ;
\draw  [draw opacity=0] (618.88,149.26) .. controls (618.96,151) and (619,152.74) .. (619,154.5) .. controls (619,218.29) and (565.18,270) .. (498.8,270) .. controls (432.41,270) and (378.6,218.29) .. (378.6,154.5) .. controls (378.6,90.71) and (432.41,39) .. (498.8,39) .. controls (561.61,39) and (613.16,85.29) .. (618.54,144.3) -- (498.8,154.5) -- cycle ; \draw   (618.88,149.26) .. controls (618.96,151) and (619,152.74) .. (619,154.5) .. controls (619,218.29) and (565.18,270) .. (498.8,270) .. controls (432.41,270) and (378.6,218.29) .. (378.6,154.5) .. controls (378.6,90.71) and (432.41,39) .. (498.8,39) .. controls (561.61,39) and (613.16,85.29) .. (618.54,144.3) ;  
\draw   (472.57,77.72) -- (443.56,51.46) -- (550,50.48) -- (525.79,77.23) -- cycle ;
\draw [color={rgb, 255:red, 74; green, 88; blue, 226 }  ,draw opacity=1 ]   (225,76.25) -- (496.01,45.71) ;
\draw [shift={(498,45.48)}, rotate = 173.57] [color={rgb, 255:red, 74; green, 88; blue, 226 }  ,draw opacity=1 ][line width=0.75]    (10.93,-3.29) .. controls (6.95,-1.4) and (3.31,-0.3) .. (0,0) .. controls (3.31,0.3) and (6.95,1.4) .. (10.93,3.29)   ;
\draw  [color={rgb, 255:red, 208; green, 2; blue, 27 }  ,draw opacity=1 ][line width=0.75] [line join = round][line cap = round] (508,41.48) .. controls (508,43.51) and (507,45.46) .. (507,47.48) ;
\draw  [color={rgb, 255:red, 208; green, 2; blue, 27 }  ,draw opacity=1 ][line width=0.75] [line join = round][line cap = round] (522,43.48) .. controls (522,44.86) and (522.37,47.48) .. (521,47.48) ;
\draw  [color={rgb, 255:red, 208; green, 2; blue, 27 }  ,draw opacity=1 ][line width=0.75] [line join = round][line cap = round] (537,45.48) .. controls (537,46.48) and (537,47.48) .. (537,48.48) ;
\draw  [color={rgb, 255:red, 208; green, 2; blue, 27 }  ,draw opacity=1 ][line width=0.75] [line join = round][line cap = round] (484,41.48) .. controls (484,44.17) and (483.1,47.58) .. (485,49.48) ;
\draw  [color={rgb, 255:red, 208; green, 2; blue, 27 }  ,draw opacity=1 ][line width=0.75] [line join = round][line cap = round] (496,42.48) .. controls (494.57,43.92) and (495,46.46) .. (495,48.48) ;
\draw  [color={rgb, 255:red, 208; green, 2; blue, 27 }  ,draw opacity=1 ][line width=0.75] [line join = round][line cap = round] (466,46.48) .. controls (466,46.82) and (466,47.15) .. (466,47.48) ;
\draw  [color={rgb, 255:red, 208; green, 2; blue, 27 }  ,draw opacity=1 ][line width=0.75] [line join = round][line cap = round] (465,45.48) .. controls (465,46.82) and (465,48.15) .. (465,49.48) ;
\draw  [color={rgb, 255:red, 208; green, 2; blue, 27 }  ,draw opacity=1 ][line width=0.75] [line join = round][line cap = round] (475,43.48) .. controls (475,45.15) and (475,46.82) .. (475,48.48) ;
\draw  [color={rgb, 255:red, 0; green, 0; blue, 0 }  ,draw opacity=1 ][line width=0.75] [line join = round][line cap = round] (619,142.25) .. controls (619,144.92) and (619,147.58) .. (619,150.25) ;
\draw  [color={rgb, 255:red, 0; green, 0; blue, 0 }  ,draw opacity=1 ][line width=0.75] [line join = round][line cap = round] (501,38.25) .. controls (500.67,38.25) and (500.33,38.25) .. (500,38.25) ;

\draw (4,148.4) node [anchor=north west][inner sep=0.75pt]    {$L_{\infty }$};
\draw (283,140.4) node [anchor=north west][inner sep=0.75pt]    {$p$};
\draw (328,41.4) node [anchor=north west][inner sep=0.75pt]    {$h_{1}$};
\draw (457,16.4) node [anchor=north west][inner sep=0.75pt]    {$q=[ 0:1:0]$};
\draw (488,90.4) node [anchor=north west][inner sep=0.75pt]    {$W_{3}$};
\draw (79,128.4) node [anchor=north west][inner sep=0.75pt]    {$V_{1}$};

\end{tikzpicture}

This set could get possibly thinner as we approach $p$ but that does not matter much for what we want. Notice this reasoning follows as $h_1$ contracts the line at infinity to $q$ and $W_3$ can be understood as a neighborhood of $q$. Analogously, there is a neighborhood $V_2$ of $L_\infty\setminus\{q\}$ with $h_1^{-1}(V_2)\subset W_4$ and neighborhoods $V_3$ and $V_4$ of $L_\infty\setminus\{\tilde{p}\}$ and of $L_\infty\setminus\{\tilde{q}\}$ with analogous properties. Set $U=\bigcap _{i=1}^4V_i$ which is open and non-empty as it contains $L_\infty\setminus\{p,q,\tilde{p},\tilde{q}\}$. Moreover, notice that $U$ was designed so that the orbit of $(x,y)\in U\cap \C^2$ goes towards infinity under any walk $S_n(\omega)=\omega_{n}\circ\dots\circ\omega_1$ such that the length of the finite word $S_n(\omega)$ goes to infinity as $n\to\infty$. More precisely: for all $(x,y)\in U$ for all $\omega\in G^\N$ with $\mbox{length}(S_n(\omega))=k$ we have

\begin{equation}\label{Eq: increasing U}
    \|S_n(\omega)\cdot(x,y)\|>2^{k-1}\|(x,y)\|.
\end{equation}

Clearly, there are sequences $\omega\in G^\N$ with $\mbox{length}(S_n(\omega))\not\to \infty$ as $n\to\infty$. However, those types of sequences are such that $(S_n(\omega))_{n\in\N}$ has a convergent subsequence (in fact a constant subsequence) by the pigeonhole principle  as we are limited by the bound over the length of $S_n(\omega)$ and we work over a finite alphabet. As normality includes also the case when the sequence of maps diverges locally to infinity, then there is always a locally uniformly convergent subsequence for any $\omega\in G^\N$ and any point in $U$, meaning $U\subset\mathcal{F}(G)$.
\end{Proof}

\begin{remark}
    Moreover, notice that if $\mbox{length}(S_n(\omega))\underset{n\to\infty}{\to} \infty$ we can go to a subsequence such that:
\begin{equation*}
    \displaystyle\lim _{k\to\infty}S_{n_k}(\omega)[x:y:1]\in\{p,q,\tilde{p},\tilde{q}\}.
\end{equation*}
This is another argument by the pigeonhole principle using that the generators contract $L_{\infty}$ only to these 4 points.
\end{remark}



     

\section{Stationary Measures}
\label{sec:2}
One way to study orbit closures is to consider more quantitative questions about random walks on a space. Breiman's Law of Large Numbers as it appears on \cite{cantatstationary} implies that the study of the distribution of the random orbits of points can be done by describing the stationary measures over the phase space $X$ induced by a measure on the group $\Gamma$ that is acting on $X$. You can also see \cite{furman} for similar remarks. With this on mind, let's define stationary measure:

\begin{definition}\cite[Chapter 8]{Einsiedler}
  Let $\Gamma$ be a $\sigma$-compact metrizable group equipped with a probability measure $\nu$. Suppose that $\Gamma$ acts continuously on a $\sigma-$compact metric space $X$. We say that $\mu\in\mathcal{M}(X)$ is $\nu$-\textbf{stationary} if  

  $$\mu=\displaystyle\int_\Gamma g_*\mu d\nu(g)$$

  or equivalently

  $$\displaystyle\int_Xf(x)d\mu(x)=\int_\Gamma\int_Xf(g\cdot x)d\mu(x)d\nu(g)$$

  for all $f\in C_c(X).$

\end{definition}

  \textbf{Special Cases:}
  
\begin{itemize}
    \item[(i)] If $\mu$ is an extreme point of the set of $\nu$-stationary measures, we call it a $\nu$-ergodic measure. 
    
    \item[(ii)]\, If for all $g\in G$ we have $g_*\mu=\mu$, we call $\mu$ a $G-$invariant measure.
\end{itemize}

We have defined stationary measure in great generality but the  basic idea is the following: Consider a group $\Gamma$ endowed with a measure $\nu$ acting on a probability space $X$ with Borel measure $\mu$. To require $\mu$ to be invariant under every element of the group is a really strong assumption, so instead we require a weaker condition: for $\mu$ to be invariant in average. Another way to understand this condition is the following: 

\begin{lemma}[{\cite[Proposition 2.2.2]{furman}}]
The measure $\mu$ is $\nu$-stationary/ergodic if and only if the measure $\nu^\N\times \mu$ on the space $\Gamma^\N\times X$ is $T$-invariant/ergodic where $T:G^\N\times X\to G^\N\times X$ is the map $T\left((\omega_n)_{n=1}^\infty, x\right):=\left(\sigma(\omega_n)_{n=1}^\infty, \omega_1\cdot x\right)$ and $\sigma$ is the left shift operator.
\end{lemma}

Our particular problem is the following: endow the group $G$ generated by H\'enon maps with the probability measure $\nu=a_1\delta_{h_1}+a_2\delta_{h_2}+a_3\delta_{h_1^{-1}}+a_4\delta_{h_2^{-1}}$ with $\sum_{i=1}^4a_i=1$  and $a_i>0$ for each $i$, i.e. a measure supported in this symmetric set of generators. Show that if $\mu$ is a $\nu-$stationary measure on $\C^2$ the support of the measure $\mu$ is compact. 

\section{Results on Stationary measures for $G$}

Now we do a transition to Random Dynamics. First we introduce an useful lemma:

\begin{lemma}
    The set of sequences $\omega\in G^\N$ such that $\left(\mbox{length}(S_n(\omega))\right)_{n\in\N}$ is unbounded has measure 1. In fact $\left(\mbox{length}(S_n(\omega))\right)_{n\in\N}$ almost surely has linear growth.
\end{lemma}
In particular this tells us that orbits of points in the set $U$ go towards infinity for almost every sequence of maps in $G$ (see equation \eqref{Eq: increasing U}). This lemma follows directly by the classical Kesten's criterion for amenability  \cite{kesten}, more precisely, the following consequence (as it appears in \cite{quint}) of this criterion:

\begin{theorem}[Kesten]\label{Theo: Kesten} Let $G$ be a discrete group and $\nu$ be a Borel probability measure on $G$. Assume that the subgroup of $G$ spanned
by the support of $\nu$ is not amenable. Let $K$ be a finite symmetric subset of $G$ with $G=\bigcup_{n\in\N}K^n$ and define a distance on $G$ by $dist(g,h):=\min\{n: g^{-1}h\in K^n\}$. Then, there exists $\epsilon>0$ such that, if $(g_n)_{n\in\N}$ is a sequence of independent random elements of $G$ with law $\nu$ we have almost surely, that for $n$ large,  $dist(\mbox{Id}, g_k\dots g_2g_1)\geq\epsilon n$.
    
\end{theorem}

In our context what this implies is that if $G=\langle h_1,h_2\rangle$ is a non-amenable group there is some $\epsilon>0$ where almost surely for $\omega\in G^\N$ we have $\mbox{length}(S_n(\omega))\geq \epsilon n$ for large $n$. This is valid for our measure $\nu$ as it is supported on a set of generators. Indeed, our group $G$ is non-amenable, moreover, $G\cong\F_2$:

\begin{lemma}
   The group $G=\langle h_1,h_2\rangle$ is isomorphic to the free group on two generators.

\end{lemma}

For this notice that the periodic behavior of the disjoint truncated cones $\{W_i\}_{1\leq i\leq 4}$ induces a Ping-pong Lemma situation (see, \cite{Loh} Chapter 4.3) as we jump between truncated cones back and forth by the appropriate map. So, $G$ is free of rank 2. Another way to see this uses the translations induced by the H\'enon maps on the Bass-Serre tree of the group of polynomial diffeomorphisms. This is slightly more technical so we refer the reader to \cite{lamythesis} for the details on Bass-Serre Theory. 

Our second main theorem is the following:

\begin{theoremA}
    If $\mu$ is a $\nu$-stationary Borel probability measure on $\C^2$ then $\supp\mu$ is compact. 
\end{theoremA}

To prove this we need the following version of the Ergodic Theorem:

\begin{theorem}[{\cite[Random Ergodic Theorem 4.1.1]{furman}}]\label{RET} Let $(X,\mu)$ be a probability space let $\Gamma$ be a group acting continuously on $X$. If $\mu$ is $\nu$-ergodic and $f\in L^1(X,\mu)$ then for $\nu^\N$-a.e. $\omega\in \Gamma^\N$ and $\mu$-a.e. $x$ we have:

\begin{equation}
  \displaystyle \lim_{N \to \infty}\frac{1}{N}\sum_{n=1}^Nf(s_n(\omega)\cdot x)=\int_X fd\mu. 
\end{equation}
    
\end{theorem}
    




We now can outline the proof of Theorem B:

\begin{Proof}
A sketch of the proof for the case of $\nu=\frac{1}{4}\left(\delta_{h_1}+\delta_{h_1^{-1}}+\delta_{h_2}+\delta_{h_2^{-1}}\right)$ reads as follows. Observe that is enough to show that the sets $U$, $W_1, W_2, W_3$ and $W_4$ have measure zero as their union is a neighborhood of $L_{\infty}$. For $U$ this follows as we can use the $\nu$-stationary property to write $\mu$ as invariant on average for reduced words of any length, in particular: 
\begin{equation}\label{eq: neat formula}
   \mu(U)=\frac{1}{\left| \Omega_n\right|}\displaystyle\sum_{\omega\in \Omega_n}\mu(\omega U) 
\end{equation}
where $\Omega_n$ is the set of reduced words of length $n$ (this because of the symmetry of the weights\footnote{It's not difficult to show, a simple induction argument works.}). Given $\epsilon>0$, we can use inner regularity to find some ball of the origin $B_{R}(0)$ such that $\mu(\C^2\setminus \overline{B}_{R}(0))<\epsilon$. Because of \eqref{Eq: increasing U} we have: for all $t>0$ exists some $k_0\in \N$ such that if $(x,y)\in U$ and $w=\omega_k\dots\omega_1\in G$ is a reduced word with $k\geq k_0$ then $w(x,y)\notin \overline{B}_t(0)$. Hence,

\begin{equation*}
    \mu(U)\leq \frac{1}{\left| \Omega_n\right|}\displaystyle\sum_{\omega\in \Omega_n}\mu(\C^2\setminus \overline{B}_{R}(0))<\frac{1}{\left| \Omega_n\right|}\displaystyle\sum_{\omega\in \Omega_n}\epsilon=\epsilon.
\end{equation*}
Thus $\mu(U)=0$ as $\epsilon>0$ was arbitrary. To show the same for the truncated cones is more delicate: Suppose that $\mu(W_1)>0$ and assume without losing generality that $\mu$ is $\nu$-ergodic (as we can go to an ergodic decomposition of $\mu$) given $\epsilon>0$ by inner regularity there exists a compact set $K_\epsilon\subset W_1$ with $0<\mu(W_1)-\epsilon<\mu(K_\epsilon)$. Then by the Random Ergodic Theorem: for $\nu^\N$-a.e. $\omega$ and $\mu$-a.e. $(x,y)$ we have

   \begin{equation}\label{eq: random}
       \displaystyle \lim_{N \to \infty}\frac{1}{N}\sum_{n=1}^N{\chi}_{K_\epsilon}(s_n(\omega)\cdot (x,y))=\mu(K_\epsilon)
   \end{equation}

 Without losing generality assume $(x,y)\in K_\epsilon$ (as this has non-zero measure) and consider sequences that start with $h_2$. This is possible as $\nu^\N(\{h_2\}\times G^\N)>0$, then \eqref{eq: random} is achieved by almost every such sequence. We highlighted that for generic $\omega$ the length of $S_n(\omega)$ grows linearly by Theorem \ref{Theo: Kesten}. This guarantees that almost surely eventually  $(x,y)$ will go out of $K_\epsilon$ and will not comeback, a contradiction. So, $\mu(W_1)=0$, an analogous reasoning works for the other cones.
 
 For the general proof: We can use Theorem \ref{RET} as we did with the truncated cones to show $\mu(U)=0$ in the general case also. Notice the reasoning for the truncated cones does not rely on the symmetry of the weights, so this works in general. 

\end{Proof}

Recall that the filled Julia set of map $h$ is $K(h)=\{(x,y):(h^n(x,y))_{n\in \Z}\,\,\mbox{bounded}\}.$ We finish by highlighting the following direct consequence:

 \begin{theoremA}
     If the filled Julia sets $K(h_1)$ and $K(h_2)$ are disjoint then there are not $\nu$-stationary measures.
 \end{theoremA}

 Indeed, suppose $(x,y)\in \supp\mu$ for some $\nu$-stationary measure. As $K(h_1)\cap K(h_2)=\emptyset$, then $(x,y)$ is not in this intersection. Without losing generality say that $(x,y)\notin K(h_1)$, meaning that there are elements of the $\Z$-orbit of $(x,y)$ under $h_1$ that can go further and further towards infinity. As $\supp\mu$ is compact (in particular bounded) then there is some $n\in\Z$ with $h_1^n(x,y)\notin\supp\mu$, then $(x,y)$ can not be in the support of $\mu$, a contradiction.
\section{Conclusion}
 
     Notice that the same ideas and results apply if we consider the group $\tilde{G}=\langle h_1,\tilde{h_2}\rangle$ where $\tilde{h_2}$ is the generic conjugated rotation of a H\'enon map $h\neq h_1$ such that $h$ is also super-attracting for $q$ and has the same indeterminacy $p$ as $h_1$. But it's expected that more sophisticated ideas are necessary for the general case. 
%
%
%

\section*{Acknowledgments}
\small\noindent
The author was partially supported by NSF grant DMS-2154414 while completing this work. I thank my advisor Dr. Roland Roeder for his constant guidance and for the interesting questions he shares with me. Also, I'm very grateful to Dr. Julio Rebelo for the useful suggestions and discussion. Finally, I also thank Arnaud Nerriere who worked independently on the same problem (and has also recently succeeded in proving the results in a more general setting) for his invaluable help and for the many interesting discussions that helped give shape to this work.


\end{document}